\documentclass[12pt,english]{article}
\usepackage[T1]{fontenc}
\usepackage[latin9]{inputenc}
\usepackage{geometry}
\usepackage{amsmath}
\usepackage{amsthm}
\usepackage{amssymb}
\usepackage{bm}
\usepackage{color}
\usepackage{todonotes}
\usepackage{graphicx}
\usepackage{bbm}
\usepackage{dsfont}
\usepackage{babel}
\usepackage{svg}
\usepackage{enumitem}

\numberwithin{equation}{section}
\numberwithin{figure}{section}
\theoremstyle{plain}
\newtheorem{thmstar}{\theoremname}
\theoremstyle{plain}
    \ifx\thechapter\undefined
      \newtheorem{prop}{\propositionname}
    \else
      \newtheorem{prop}{\propositionname}[chapter]
    \fi
\theoremstyle{plain}
    \ifx\thechapter\undefined
      \newtheorem{lem}{\lemmaname}
    \else
      \newtheorem{lem}{\lemmaname}[chapter]
    \fi
\theoremstyle{remark}
\newtheorem{rem*}{\remarkname}
\newtheoremstyle{named}{}{}{\itshape}{}{\bfseries}{.}{.5em}{\thmnote{#3}#1}
\theoremstyle{named}
\newtheorem{thm}{Proposition 1'}
\newenvironment{thmbis}[1]
{\renewcommand{\thethm}{\ref{#1}$'$}%
  \addtocounter{thm}{-1}%
  \begin{thm}}
  {\end{thm}}

\newtheorem{rem}{Remark}
\newtheorem{cor}{Corollary}

\providecommand{\lemmaname}{Lemma}
\providecommand{\propositionname}{Proposition}
\providecommand{\remarkname}{Remark}
\providecommand{\theoremname}{Theorem}

\newcommand{\N}{\mathbb{N}}
\newcommand{\R}{\mathbb{R}}
\newcommand{\Z}{\mathbb{Z}}
\newcommand{\E}{\mathbb{E}}
\newcommand{\F}{\mathbb{F}}

\newcommand{\bl}{\boldsymbol{l}}
\newcommand{\bL}{\bm{L}}
\newcommand{\bu}{\boldsymbol{u}}
\newcommand{\bv}{\boldsymbol{v}}
\newcommand{\bp}{\bm{p}}
\newcommand{\bx}{\boldsymbol{x}}
\newcommand{\by}{\boldsymbol{y}}
\newcommand{\bz}{\boldsymbol{z}}

\newcommand{\fqed}{\tag*{\ensuremath{\square}}}
\renewcommand{\qed}{\hfill\square}
\newcommand{\upinf}{\uparrow \infty}
\newcommand{\iid}{\overset{\text{iid}}{\sim}}

\renewcommand{\P}{\mathbb{P}}
\renewcommand{\L}{\mathbb{L}}

\newcommand{\eL}{\mathcal{L}}
\newcommand{\frC}{\mathfrak{C}}
\newcommand{\frL}{\mathfrak{L}}
\newcommand{\eG}{\mathcal{G}}

\newcommand{\eF}{\mathcal{F}}
\newcommand{\eB}{\mathcal{B}}
\newcommand{\Cset}{\mathcal{C}}
\newcommand{\Cyl}[1]{\text{Cyl}(#1)}

\newcommand{\eLind}[2]{\mathcal{L}_{#1}^{#2}}

\newcommand{\floo}[1]{\lfloor#1\rfloor}
\newcommand{\ceil}[1]{\lceil#1\rceil}
\newcommand{\ip}[1]{\langle#1\rangle}

\newcommand{\ind}{\mathds{1}}

\usepackage[affil-it]{authblk}

\title{Shape theorem and surface fluctuation for Poisson cylinders}
\author[1]{Marcelo Hilario}
\author[2]{Xinyi Li}
\author[2]{Petr Panov}
\affil[1]{\small Universidade Federal de Minas Gerais}
\affil[2]{\small The University of Chicago}

\date{\small \today}

\makeatletter
\renewcommand{\@biblabel}[1]{[#1]\hfill}
\makeatother

\newcommand{\cgrid}{c_0}
\newcommand{\cdist}{c_1}
\newcommand{\cpatch}{c_2}
\newcommand{\clocal}{c_3}
\newcommand{\cgrow}{c_4}
\begin{document}

\maketitle
\begin{abstract}
  In this work, we prove a shape theorem for Poisson cylinders and give a power law bound on surface fluctuations. We prove that for any $a \in (1/2, 1)$, conditioned on the origin being in the set of cylinders, every point in this set, whose Euclidean norm is less than $R$, lies at an internal distance less than $R + O(R^a)$ from the origin.
\end{abstract}

\section{Introduction}\label{sec:intro}
We consider a random collection of bi-infinite cylinders, which are sampled by thickening the elements in the support of a Poisson point process in the space of all the lines in $\R^d$, $d\geq 3$. This model, called the \emph{Poisson cylinder model} serves as a natural mathematical model for various random fiber structures and has many applications in image analysis. We refer the readers to \cite{Spiess} for a detailed survey on this topic. Recently, many  geometric properties of this model, especially those related to percolation, have been studied \cite{tyke_win, hilario15, tyke_bro} along with other models presenting long-range correlation.

We now describe our main result. Fix $d \geq 3$ and let $\P^u$ stand for the law of Poisson cylinders in $\R^d$ with intensity measure $u\mu$ where $u>0$ and  $\mu$ is the translation and rotation invariant Haar measure on the space of lines in $\R^d$. We refer the readers to Section \ref{sec:model} for a precise mathematical construction. Let $\Cset \subset \R^d$ be the union of the cylinders and $\rho = \rho(\Cset)$ the random metric defined as follows: for two points $a, b \in \Cset$, $\rho(a, b)$ equals the minimal length of  paths connecting $a$ and $b$ that stay entirely inside $\Cset$. We call this metric on $\Cset$ the \emph{internal distance}.

Let $B_r$ and $B_r^\rho$ stand for the closed ball of radius $r$ centered at the origin $\bm 0$  with regard to the Euclidean metric and $\rho$, respectively. 
Our main result is
\begin{thmstar}[Shape theorem]
  For every $a \in (1/2, 1)$, $u>0$ and $c > 0$, $\P^u[\cdot|\bm 0 \in \Cset]$-almost surely there exists a finite $R_0 > 0$ such that 
  \begin{equation} \label{eq:shapethm}
    (\Cset \cap B_R) \subseteq B_{R + c R^a}^\rho \ \ \mbox{ for all } \ \ R > R_0.
  \end{equation}
\end{thmstar}
This theorem can be roughly rephrased as follows: given a realization of the set of cylinders such that $\bm 0 \in \Cset$, for large $R$ every point in $B_R \cap \Cset$ can be connected to $\bm 0$ by a path that lies inside $\Cset$, and whose total length is bounded from above by $R\bigl( 1 + O(R^{1-a}) \bigr)$. Note that conditioning on $\bm 0 \in \Cset$ is an arbitrary choice: by the translational invariance of the cylinder model, we could have taken any other point in $\R^d$.
\begin{rem}
  Our result is stronger than usual shape theorems, which, in the notation of this work, would be stated as follows: there exists a convex compact $D \subset \R^d$, called the asymptotic shape with respect to $\rho$, such that for any $\epsilon > 0$,   $\P[\cdot|\bm 0 \in \Cset]$-almost surely, there is a finite $R_0 > 0$ such that 
  \begin{equation} \label{eq:usualshapethm}
    B_{(1-\epsilon)R}^\rho \subseteq ( \Cset \cap D_R ) \subseteq B_{(1+\epsilon)R}^\rho \ \mbox{ for } R\geq R_0,
  \end{equation}
  where 
  \[
    D_R = \{ a \bx \in \R^d;\, a \in [0, R], \, \bx \in D \}.
  \]
  In other words, not only do we prove that the asymptotic shape is the unit ball, but also the asymptotic equivalence between $\rho$ and the Euclidean metric, providing a bound on surface fluctuations which is similar to results in {\rm \cite{Kesten}}. For further discussion on shape theorems for models of first passage percolation, see Section 3 of {\rm \cite{50years}}. Note also that the first inclusion in \eqref{eq:usualshapethm} is immediate in our case. 
\end{rem}
We do not discuss the Poisson cylinder model for $d=2$. In fact, on a plane it suffices to look at the collection of lines, as thickening lines into cylinders is no longer necessary to guarantee connectedness. Geodesics in models of this type were first investigated by Aldous and Kendall, who proved results that amount to an $O(\log R)$ surface fluctuation in the shape theorem. See e.g.\ \cite{AK} and \cite{Kendall} for more details. It is also noteworthy that the most natural generalization of their model for higher dimensions is not the Poisson cylinder model, but rather the Poisson flats model (introduced in \cite{flats}), which is a Poissonian soup of $(d - 1)$-dimensional affine spaces. In this case the fluctuation in the shape theorem is also of order $O(\log R)$ by a projection argument.

One can also compare our result with shape theorems obtained for discrete percolation models with long-range correlations, e.g., random interlacements and level sets of the Gaussian free field. In \cite{chem_2}, a common scheme is developed for proving shape theorems for these models, which involves checking that the specific model under consideration fulfills a few criteria. Once these criteria are met, one has a shape theorem (in the form of the one in the remark above) for this model, along with lots of other geometric properties. However, it cannot be applied to Poisson cylinders due to the spatial rigidity of cylinders: for (a discretized version of) the Poisson cylinder model, assumption {\bf P3} of \cite{chem_2} which is known as the decoupling inequality in the random interlacement folklore, is not satisfied. Since our result in this work is actually stronger than statement in the remark above, we will not seek a sophisticated adaptation of \cite{chem_2} to bypass this obstacle.

We now explain the strategy of the proof of \eqref{eq:shapethm}. Precise statements and detailed explanations can be found in Section \ref{sec:res}. 

$\bullet$ As a first step, we reduce the original theorem to a statement regarding the internal distance between $\bm 0$ and a point $\bx \in \Cset \cap B_R$. This is summarized in Proposition \ref{p:bound_rho}.

$\bullet$ Then, we show that the $\mu$-measure of a ``\emph{local network}'' of truncated cylinders with length of order $r = R^a$ near $\bm 0$ (resp.\ $\bx$) is in some sense comparable to that of a Euclidean ball of radius $r$. See Proposition \ref{p:local} for a precise statement.

$\bullet$ Finally, we find a ``\emph{highway}'' (long cylinder) connecting the local networks near $\bm 0$ and $\bx$. Thanks to the previous step we know that local networks are about  as "visible" as Euclidean balls of the same size which, together with a classical estimate on $\mu$ (see Lemma 3.1 in \cite{tyke_win}) assures the existence of a highway with high probability. This part corresponds to Lemma \ref{l:global}.

One can compare our result with the connectivity results of Poisson cylinders. Imagine a graph where each vertex represents a cylinder in the Poissonian soup, and where edges connect any two intersecting cylinders. In \cite{RIconnectivity} the authors show that for any intensity $u>0$ this graph is $\P^u$-almost surely well-connected and its diameter is equal to $(d-1)$. However, their results do not provide a bound on $\rho$. On the other hand, our strategy, which also involves creating connections between cylinders, provides a short path, but it is not intended to optimize the amount of cylinders visited by this path.

It is worth mentioning that the second and third steps above lead to a strong connectivity result in the form of criterion {\bf S1} in \cite{chem_2} for Poisson cylinders. As a comparison, see \cite{transience} and, in particular, Lemma 12 therein, for corresponding results for random interlacements. In their case, $\mu$ and ``visibility'' for cylinders are replaced by random walk capacity and hitting probability, respectively.

With our proof strategy, the lower bound for $a$ cannot be improved further. Indeed, if $a \leq 1/2$, then the local networks from the construction above will no longer be visible to each other with high probability. However, we are not able to rule out the possibility of a completely different strategy which could lead to stronger results. For instance, it is not impossible that the shortest paths inside the set of cylinders consist not of a single long cylinder and a few short ones, but rather of many shorter segments, much like as in the models investigated by Aldous and Kendall.

This work allows for various extensions. If $\sqrt{R}$ is indeed the right scaling for surface fluctuations, we are naturally led to the question of whether more can be said on these fluctuation. To mix this problem with classical Bernoulli first passage percolation problems, one can assign random speed on each cylinder, or even between different sections of the same cylinder, and ask a similar question. One can also ask if the same shape theorem holds for the Poisson cylinder set in hyperbolic space, where the connectedness of cylinders undergoes a phase transition as the soup intensity changes, see \cite{tyke_bro_hyper} for more details.

We now explain how this work is organized. In Section \ref{sec:model} we introduce the model and our notations. In Section \ref{sec:res} we state our main result and a few key propositions. Proofs are postponed till Section \ref{sec:proof}.

{\bf Acknowledgments:} Part of this work was accomplished during a visit of XL to NYU Shanghai where MH was a long-term visitor. MH and XL would like thank Vladas Sidoravicius for warm hospitality and useful discussions. The authors would like to thank Antonio Auffinger for useful comments on the text. XL would also like to thank Yuval Peres for pointing out a reference. MH was supported by CNPq grants ``Projeto Universal" (307880/2017-6) and ``Produtividade em Pesquisa" (406659/2016-8) and by FAPEMIG grant ``Projeto Universal" (APQ-02971-17). 

\section{Model, notation and conventions}\label{sec:model}
In this section, we introduce notations and describe the Poisson cylinder model.

\subsection{Notation}
Throughout this work, we consider $\R^d$ with $d \geq 3$ and view the integer lattice $\Z^d$ as its subset. Ordered tuples and particularly vectors in $\R^d$ are written in bold. We use $|\cdot|$ for the Euclidean norm on $\R^d$. We denote by $\bx + B_r$ the closed Euclidean ball of radius $r>0$ centered at $\bx \in \R^d$; here the plus sign stands for the sumset operation and $\{ \bx \}$ is replaced by $\bx$ for brevity. 
Given a metric $\rho: \R^d \times \R^d \mapsto [0, \infty]$, we let 
\[
  B_r^\rho (\bx) = \bigl\{ \by \in \R^d: \rho(\bx, \by)\leq r \bigr\} \subset \R^d,
\]
and write $B_r^{\rho} =B_r^{\rho}(\bm 0)$ for simplicity.

We write $(\cdot_j)_{j \in A}$ for a sequence of elements whose indices take values in an ordered countable set $A$, and $\{ \cdot_j \}_{j \in A}$ for unordered sets with any index set $A$. If a set $A$ is finite, we use $|A|$ for its cardinality. For $x \in \R$, we write
\[
  \floo x = \sup\{y \in \Z: y \leq x\}, \qquad \ceil x = \inf\{y \in \Z: y \geq x\},
\]
and use $[x]$ to denote the set $\{1, 2, \ldots, \floo x\}$.

We will use a number of positive and finite constants which will be denoted by $c$ and whose values might change from line to line. Even when we do not write it explicitly, these constants will always assume strictly positive values. When $c$ comes with an integer subscript, its value is kept fixed throughout the paper. Symbolic superscripts are exponents and not indices. For example, $\cgrow^M$ refers to a fixed constant $\cgrow>0$ raised to the power $M$.

Given two transformations $f$ and $g$ of $(0, \infty)$, we write $f \in O(g)$ and $g \in \Omega(f)$, if there is a $c$ such that $f \leq c g$. If $f \in O(g)$ and $g \in O(g)$, then we write $f \asymp g$.

Finally, unless otherwise specified, $\log$ stands for the natural logarithm.

\subsection{Poisson cylinder model}
We now turn to the model of Poisson cylinders. Let $\L$ be the set of 1-dimensional affine subspaces of $\R^d$. Fix any line $\hat{l} \in \L$. For every line $l \in \L$ there exist a translation $\tau$ by a vector orthogonal to $\hat{l}$ and a rigid rotation $\theta$ around the origin, such that $(\theta \circ \tau)(\hat{l}) = l$. Endow $\L$ with the finest topology which makes $\theta \circ \tau$ continuous for all $\tau$ and $\theta$. Once the topology is given, we equip $\L$ with the Borel $\sigma$-algebra $\eB(\L)$ and denote by $\mu$ the Haar measure on $\bigl(\L, \eB(\L) \bigr)$ which is invariant under all isometries of $\R^d$. The measure $\mu$ is unique up to a normalizing constant, but we do not choose the latter explicitly, because our results do not depend on the normalization.

For $l \in \L$, define a cylinder of radius $1$ around $l$ as $\Cyl{l} = l + B_1$. For open or compact $A \subset \R^d$ let 
\[
  \eL(A) = \{l \in \L: \Cyl{l} \cap A \neq \emptyset\},
\]
which is $\eB(\L)$-measurable; see \cite{tyke_win}, around Eq.\ (2.11) for a proof.

Another important Borel measurable subset of $\L$ is the set of lines whose angle with respect to a given vector falls within a certain range. For $\bx \in \R^d$, a unit vector $\bu \in \R^d$ and $0 \leq \alpha \leq \beta \leq \pi/2$, we let
\[
  \eLind{\bx, \bu}{\alpha, \beta} = \bigl\{ l \in \eL \{ \bx \}: |\ip{\bu, \bv(l)}| \in [\sin \alpha, \sin \beta] \bigr\},
\]
where $\bv(l)$ is either of the two unit vectors directing a line $l \in \L$, and $\ip{\cdot, \cdot}$ is the scalar product. For example, $\eLind{\bm 0, \bv(l_0)}{0, 0}$ is the set of lines hitting $B_1$ and orthogonal to $l_0$.

We now consider the Poisson point process on $\bigl(\L, \eB(\L) \bigr)$ defined on a probability space $(\mathcal{M}, \mathcal{A}, \P^u)$. Here
\[
  \mathcal{M} = \Big\{ \omega = \sum_{j \geq 0} \delta_{l_j} : \, l_j \in \L \ \text{ for all $\ j,\ $ and } \ \omega(A) < \infty \ \text{ for compact } \ A \in \eB(\L)\Big\}
\]
is the sample space (composed of the locally finite point measures). The set of events $ \mathcal{A} = \sigma \bigl( \{e_A \}_{A \in \eB(\L)} \bigr)$ is the $\sigma$-algebra generated by the evaluation maps $e_A: \omega \mapsto \omega(A)$. Finally, $\P^u$ is the probability measure under which $\omega$ is a Poisson point process on $\L$ with intensity measure $u\mu$ for some $u > 0$. Note that $\P^u$ inherits the invariance under translations and rotations from $\mu$. See \cite{tyke_win} for a more detailed account of the properties of $\P^u$.

As mentioned before, the statements in this paper hold for any intensity parameter $u > 0$ and dimension $d \geq 3$, so we will often not be explicit about the dependence of constants and probability measures on them. For example, we write $\P(\text{event}) < c$, if for any $u>0$ there is a constant $c=c(u, d)>0$ such that $\P^u(\text{event})<c$. The same convention applies to the asymptotic notations that we have introduced previously.

Having constructed the Poisson point process of lines, we denote the set of cylinders by
\[
  \Cset = \Cset(\omega) = \bigcup_{l \in \omega}\Cyl{l}.
\]
Here and in what follows we write $l \in \omega$ and $\Cyl{l} \in \omega$ instead of $l \in \text{supp}(\omega)$ for brevity. Note that $\Cset$ is invariant in law under isometries of $\R^d$.

Given some $\bx, \by \in \R^d$, we denote by $[\bx, \by]$ the line segment connecting $\bx$ and $\by$:
\[
  [\bx, \by] = \big\{ (1-t)\, \bx + t\, \by: t \in [0, 1] \big\}.
\]
For any $A \subseteq \R^d$, the set of \emph{polygonal paths} from $\bx$ to $\by$ in $A$ is a set of finite sequences of vertices, for which the line segments between consecutive elements are within $A$:
\[
  \mathcal{P}_A(\bx, \by) = \bigl\{ (\bz_j)_{j=0}^n: n \geq 1,\ \bz_0 = \bx,\ \bz_n = \by \ \text{ and } \ [\bz_{k-1}, \bz_k] \subseteq A \quad \forall k \in [n] \bigr\}.
\]
We define the \emph{internal distance} $\rho = \rho (\omega)$ as follows:
\[
  \rho(\bx, \by) = \inf \left\{ \sum_{j=1}^n |\bz_j - \bz_{j-1}|: (\bz_j)_{j=0}^n \in \mathcal{P}_\Cset(\bx, \by) \right\}
\]
for all $\bx, \by \in \R^d$. We follow the convention that $\inf\{\emptyset\} = +\infty$. Sometimes we write $\rho = \rho(\omega)$ or $\rho = \rho(\mathcal{C})$ in order to indicate the dependence on the underlying Poisson process. It is proved in Theorem~6.1 from \cite{tyke_bro}, that $\Cset(\omega)$ is a connected set for $\P$-almost all $\omega$. In other words, $\rho (\bx, \by) < \infty \iff \bx, \by \in \Cset$ for $\bx, \by \in \R^d$.
\section{Proof strategy and intermediate results}\label{sec:res}
In this section, we structure the proof of the shape theorem and state intermediate results. First, in Section \ref{sec:3.1} we reduce our theorem to the study of the internal distance between a pair of points. Then we divide the problem into finding a ``highway'' and ``local networks'' in Section \ref{sec:3.2}. Section \ref{sec:3.3} is dedicated to the study of local networks.

We start by restating \eqref{eq:shapethm} in a form which is easier to analyze. 
\begin{thmstar}[Shape theorem restated]
  For every $a \in (1/2, 1)$ and $c>0$, $\P[\cdot|\bm 0 \in \Cset]$-a.s.
  \begin{equation*}
    \exists\, R_0 > 0 \ \text{ such that } \ \rho(\bm 0, \bx) \leq R + c R^a, \qquad \forall\, \bx \in \Cset\cap B_R, \quad \forall\, R > R_0.
  \end{equation*}
\end{thmstar}

\subsection{A bound on the internal distance} \label{sec:3.1}
To prove \eqref{eq:shapethm}, we first show that with high probability the internal distance between any pair of points $\Cset$ cannot be much larger than the Euclidean distance between them.
\begin{prop} \label{p:bound_rho}
  Given $\delta \in (0, 1)$ such that $\delta < (2a - 1)(d - 1)$, there exists $\cdist$ such that,
  \begin{equation} \label{eq:bound_rho1}
    \P \bigl[ \rho(\bm 0, \bx) > |\bx| + \cdist|\bx|^a \bigl| \bm 0, \bx \in \Cset \bigr] \in    O \bigl( \exp \{ -|\bx|^\delta \} \bigr)
  \end{equation}
  for all
  \begin{equation} \label{eq:xrange}
    \bx \in \frac{1}{2\sqrt{d}} (\R^d \backslash B_1).
  \end{equation}
\end{prop}
\begin{rem}
  Proposition~\ref{p:bound_rho} is stronger than what we really need. It suffices to provide a bound like $O(|\bx|^{-d - 1 - c})$ on the right-hand side of \eqref{eq:bound_rho1} with some  $c>0$.
\end{rem}
Pick an $\bx$ that satisfies \eqref{eq:xrange} and write $R = |\bx|$ and $r = R^a$. To avoid conditioning on the event $\{ \bm 0, \bx \in \Cset \}$, we prove an unconditional version of Proposition~\ref{p:bound_rho}, which gives a uniform bound over all possible pairs of  $l_0 \in \eL \{ \bm 0 \}$ and $l_x \in \eL \{ \bx \}$. To state it precisely, given such $l_0$, $l_x$ and any $c>0$ let
\[
  \omega^- = \omega - \omega\, \ind_{\eL \{ \bm 0, \bx \}}, \qquad \Cset_{l_0, l_x}(\omega) = \Cset(\omega^-) \cup \Cyl{l_0} \cup \Cyl{l_x},
\] 
and define the event
\begin{equation} \label{e:def_E}
  E_{l_0, l_x}(c) = \left\{ \forall\, (\bz_j)_{j=0}^n \subset \mathcal{P}_{\Cset_{l_0, l_x}(\omega)}(\bm 0, \bx),\ \ \sum_{j=1}^n|\bz_j - \bz_{j-1}| > R + c r \right\}.
\end{equation}
\begin{thmbis}{p:primed}
  For any $\delta$ as in Proposition \ref{p:bound_rho} there is $\cdist>0$ such that  \begin{equation} \label{e:bound_rho2}
    -\log \sup_{l_0 \in \eL \{ \bm 0 \},\, l_x \in \eL \{ \bx \}} \P \bigl( E_{l_0, l_x} (\cdist) \bigr) \in O(R^\delta).
  \end{equation}
\end{thmbis}

\subsection{Highway and local connections} \label{sec:3.2}

Although Propositions \ref{p:bound_rho} and \ref{p:bound_rho}' do not specify how to construct a short path, the bound on $\rho$ strongly suggests the following strategy: \begin{itemize}[noitemsep, topsep=0pt]
  \item take $B_r$ and $B_r + \bx$;
  \item if a cylinder intersects both balls, we call it a \emph{highway};
  \item if a highway exists, and we manage to connect $\bm 0$ and $\bx$ to it by finitely many truncated cylinders with heights of order $r$, then $\rho(\bm 0, \bx)$ is bounded above by $R + \cdist r$ for some $\cdist>0$.
\end{itemize}
It is not difficult to see that a highway exists. Thanks to Lemma 3.1 in \cite{tyke_win}, 
\[ 
  \mu \bigl( \eL(B_r) \cap \eL(B_r + \bx) \bigr) \in \Omega \bigl( (r^2/R)^{d-1} \bigr).
\]
Since $a > 1/2$, the probability of having no highway  between $B_r$ and $B_r + \bx$ decays very fast as we increase $|\bx|$.

The purpose of this subsection is to set up a structure in space (referred to in the text as a {\it crossing}) which mediates connections between local networks and the highway. We will often ``discretize'' various geometric objects in $\R^d$ like disks or cylinders, by replacing them with finite sets of points. Doing so allows us to use finite sequences of random variables to provide the concentration bounds from the previous subsection.

\begin{figure}[ht!]
  \centering
  \def\svgwidth{\linewidth}
\begingroup%
  \makeatletter%
  \providecommand\color[2][]{%
    \errmessage{(Inkscape) Color is used for the text in Inkscape, but the package 'color.sty' is not loaded}%
    \renewcommand\color[2][]{}%
  }%
  \providecommand\transparent[1]{%
    \errmessage{(Inkscape) Transparency is used (non-zero) for the text in Inkscape, but the package 'transparent.sty' is not loaded}%
    \renewcommand\transparent[1]{}%
  }%
  \providecommand\rotatebox[2]{#2}%
  \ifx\svgwidth\undefined%
    \setlength{\unitlength}{493.54215396bp}%
    \ifx\svgscale\undefined%
      \relax%
    \else%
      \setlength{\unitlength}{\unitlength * \real{\svgscale}}%
    \fi%
  \else%
    \setlength{\unitlength}{\svgwidth}%
  \fi%
  \global\let\svgwidth\undefined%
  \global\let\svgscale\undefined%
  \makeatother%
  \begin{picture}(1,0.392)%
    \put(0,0){\includegraphics[width=\unitlength,page=1]{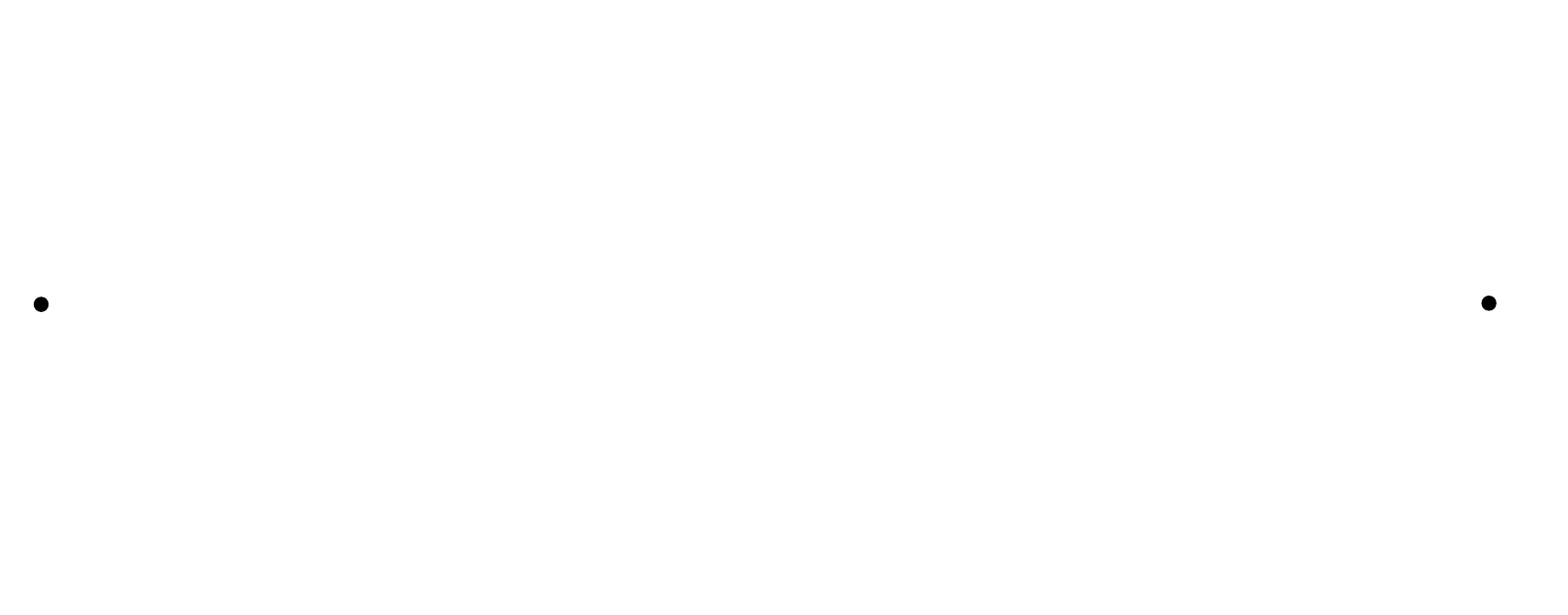}}%
    \put(0.026,0.176){\color[rgb]{0,0,0}\makebox(0,0)[lb]{$\bm 0$}}%
    \put(0.944,0.209){\color[rgb]{0,0,0}\makebox(0,0)[lb]{$\bx$}}%
    \put(0,0){\includegraphics[width=\unitlength,page=2]{general.pdf}}%
    \put(0.223,0.256){\color[rgb]{0,0,0}\makebox(0,0)[lb]{$\cgrid$}}%
    \put(0,0){\includegraphics[width=\unitlength,page=3]{general.pdf}}%
    \put(0.855,0.223){\color[rgb]{0,0,0}\makebox(0,0)[lb]{$\cpatch r$}}%
    \put(0,0){\includegraphics[width=\unitlength,page=4]{general.pdf}}%
    \put(0.015,0.116){\color[rgb]{0,0,0}\makebox(0,0)[lb]{$l_0$}}%
    \put(0,0){\includegraphics[width=\unitlength,page=5]{general.pdf}}%
    \put(0.7,0.36){\color[rgb]{0,0,0}\makebox(0,0)[lb]{$l_x$}}%
    \put(0,0){\includegraphics[width=\unitlength,page=6]{general.pdf}}%
    \put(0.219,0.369){\color[rgb]{0,0,0}\makebox(0,0)[lb]{$\frC_0$}}%
    \put(0.838,0.056){\color[rgb]{0,0,0}\makebox(0,0)[lb]{$\frC_x$}}%
    \put(0,0){\includegraphics[width=\unitlength,page=7]{general.pdf}}%
    \put(0.758,0.09){\color[rgb]{0,0,0}\makebox(0,0)[lb]{$r$}}%
    \put(0,0){\includegraphics[width=\unitlength,page=8]{general.pdf}}%
  \end{picture}%
\endgroup%
  \caption{Dashed lines represent lines in the soup $\omega$ (more specifically, here we represent the two ``given'' lines $l_0$ and $l_x$, one line from the local network near $\bm 0$ and the highway). The solid broken line depicts a polygonal path from $\mathcal{P}_{\Cset}(\bm 0, \bx)$.  The unions of filled and empty dots on disks are the crossings $\frC_0$ and $\frC_x$. The filled dots correspond to $\frC_0'$ and $\frC_x'$, and they are used to connect a highway.}
  \label{fig:general}
\end{figure}

For some big $\cpatch > 0$, which will be picked in Lemma \ref{l:last}, put two identical $(d-1)$-dimensional disks of radius $r$ between $\bm 0$ and $\bx$, at a distance $\cpatch r$ from them. We put a square grid of points on each one of those disks with mesh  $\cgrid := 10d$. We call these grids \emph{crossings}. More precisely, we take an orthogonal coordinate system so that $\bm x$ points in the direction of the last basis vector, and define the crossings as follows:
\[
  \frC_0 = \left\{ \bigl( \cgrid\, \Z^{d - 1} \cap B_r \bigr) \times \{ 0 \} + \cpatch r \, \frac{\bx}{R} \right\}, \ \ \frC_x = \left\{ \bigl( \cgrid\, \Z^{d - 1} \cap B_r \bigr) \times \{ 0 \} + \Bigl(1 - \cpatch r \, \frac{\bx}{R} \Bigr) \right\},
\]
see Fig. \ref{fig:general} for a sketch.

Note that $|\frC_j| \asymp r^{d-1}$ for $j \in \{0,x\}$. The next statement implies that under certain conditions a highway exists with high probability.
\begin{lem} \label{l:global}
  Suppose that $\frC_j' \subseteq \frC_j$, $j\in\{0,x\}$ satisfy $|\frC_j'| \in \Omega(r^\chi)$ with $\chi(\epsilon) := (d - 1) (1 - \epsilon)$ for some  $\epsilon > 0$ such that $2a (1 - \epsilon) > 1$. Then 
  \[
    -\log \P \bigl\{ \omega^- \bigl( \eL(\frC_0') \cap \eL(\frC_x') \bigr) = 0 \bigr\} \in \Omega (r^{2\chi}/R^{d-1}).
  \]
\end{lem}
Note that for $d\geq3$ our restriction on the values of $a$ and $\epsilon$ implies that $\chi(\epsilon)>1$. It now suffices to show that $l_0$ can be connected to at least $r^\chi$ points on $\frC_0$ via finitely many truncated cylinders of length $O(r)$ with high probability. A similar statement about a local network near $\bx$ will then follow, thanks to the translation and rotation invariance of the model.
\begin{prop} \label{p:local}
  For any $\epsilon$ and $\chi(\epsilon)$ be given as in Lemma \ref{l:global} there exists $\clocal = \clocal(\epsilon)$ such that for any $l_0 \in \eL \{ \bm 0 \}$,
  \begin{equation} \label{e:local}
    -\log \P \left\{ |\frC_0 \cap B_{\clocal r}^{\rho(\omega^-+\delta_{l_0})}| < r^\chi \right\} \in \Omega(r).
  \end{equation}
\end{prop}
Proposition \ref{p:local} ensures, that with high probability we have sufficiently many points to hook onto in Lemma \ref{l:global}, and those points are not too far from $\bm 0$ in $\Cset$. It implies that the $\mu$-measure of a local network of size $r$ is in $\Omega(r^\chi)$, where $\chi$ can be made arbitrarily close to $d-1$. Therefore, the "visibility" of a local network is asymptotically close to that of a Euclidean ball. The purpose of the next subsection is to outline the proof of this.

\subsection{Local networks} \label{sec:3.3}
Here we describe the strategy used for the proof of \eqref{e:local}.

We first fix some positive integer $M$. In order to exploit some independence, we split the original cylinder process into a union of $(M+1)$ independent cylinder processes with intensity $w = u/(M+1)$. Next we draw concentric spheres having radius of order $r$. We then sample cylinders from the first portion to connect $l_0$ to the first sphere and then sample cylinders from the second portion to connect the first sphere to the second one. Here, \emph{connecting} a set $A$ to a set $B$ means that there is at least one cylinder $C \in \omega^-$ such that $A \cap C \neq \emptyset$ and there exists a point in $B \cap C$, which is referred to as a \emph{connection}. We proceed iteratively, until the $M^\text{th}$ sphere is reached. We guarantee that the number of connections grows rapidly between each step with high probability. In particular, after the $M^\text{th}$ step we have sufficiently many points near the last sphere which are connected to $\Cyl{l_0}$ by a path of length in $O(r)$, to ensure that the last portion contains enough cylinders to connect those points to at least $r^\chi$ points on $\frC_0$.

More precisely, fix $\epsilon>0$ such that $2a(1-\epsilon)>1$ as in Lemma \ref{l:global} and let
\[
  M = \ceil{\log_{\frac{d - 2}{d - 1}}\epsilon} + 1, \quad a_j = \Big( \frac{d - 2}{d - 1} \Big)^j, \quad b_j = (d - 1)(1 - a_j), \qquad \forall j\in[M].
\]
Here $M$ corresponds to the number of steps as described above. The order of the length of the truncated cylinders (which are close to the $j^\text{th}$ sphere) used to form the connection on the next step is $r^{a_j}$. The order of number of connections to $\Cyl{l_0}$ provided at step $j$ is $r^{b_j}$. Note that $M$ is chosen so that $a_M < \epsilon$ and $b_M > \chi$. Also, $b_{j+1} = a_j + b_j$ for all $j \in [M-1]$.

We fix a constant $\cgrow$, whose actual value will be given by Lemma \ref{l:intermediate} below. Define the annuli $A_j := B_{\cgrow^jr} \backslash B_{\cgrow^jr-r^{a_j}}$ for $j\in[M]$. Note that if $\cgrow>10$, then for any $l\in\eL(B_r)$ the length of both segments of $l \cap B_{\cgrow r} \backslash B_{(\cgrow-1)r}$ is bounded above by $2r$ for every $r$ large enough. We will assume that as we move on. 

A sequence of points on a line segment $[\by, \bz] \subset \R^d$ is called a $\cgrid$\emph{-grid}, if the spacing between consecutive points along the segment is equal to $\cgrid$ and, in addition, no more points can be added to the sequence without violating this spacing constraint. We take a $\cgrid$-grid on $l_0 \cap B_r$ and exclude $\bm 0$ from it:
\[
  (\bl_{0,1}^m)_{m=1}^{N_{0,1}} = \cgrid \bigl( \Z \cap [-r/\cgrid, r/\cgrid] \backslash \{\bm 0\} \bigr);
\] 
note that we write $N_{0,1}$ for the amount of points in the sequence. We let $\bL_0 = \bigl( (\bl_{0,1}^m)_{m=1}^{N_{0,1}},l_{0,1} \bigr)$ and call it the \emph{zeroth layer}.

We then recursively stack a collection of geometric structures (essentially, sets of points) onto the zeroth layer, and we also call them \emph{layers}. Layers consist of ``beaded'' \emph{threads}, each thread being a pair formed by a $\cgrid$-grid and the line on which the grid lies.

\begin{figure}[h!]
  \begin{minipage}[c]{.5\textwidth}
    \centering
    \def\svgwidth{1\textwidth}
\begingroup%
  \makeatletter%
  \providecommand\color[2][]{%
    \errmessage{(Inkscape) Color is used for the text in Inkscape, but the package 'color.sty' is not loaded}%
    \renewcommand\color[2][]{}%
  }%
  \providecommand\transparent[1]{%
    \errmessage{(Inkscape) Transparency is used (non-zero) for the text in Inkscape, but the package 'transparent.sty' is not loaded}%
    \renewcommand\transparent[1]{}%
  }%
  \providecommand\rotatebox[2]{#2}%
  \ifx\svgwidth\undefined%
    \setlength{\unitlength}{321.52721036bp}%
    \ifx\svgscale\undefined%
      \relax%
    \else%
      \setlength{\unitlength}{\unitlength * \real{\svgscale}}%
    \fi%
  \else%
    \setlength{\unitlength}{\svgwidth}%
  \fi%
  \global\let\svgwidth\undefined%
  \global\let\svgscale\undefined%
  \makeatother%
  \begin{picture}(1,1.00004261)%
    \put(0,0){\includegraphics[width=\unitlength,page=1]{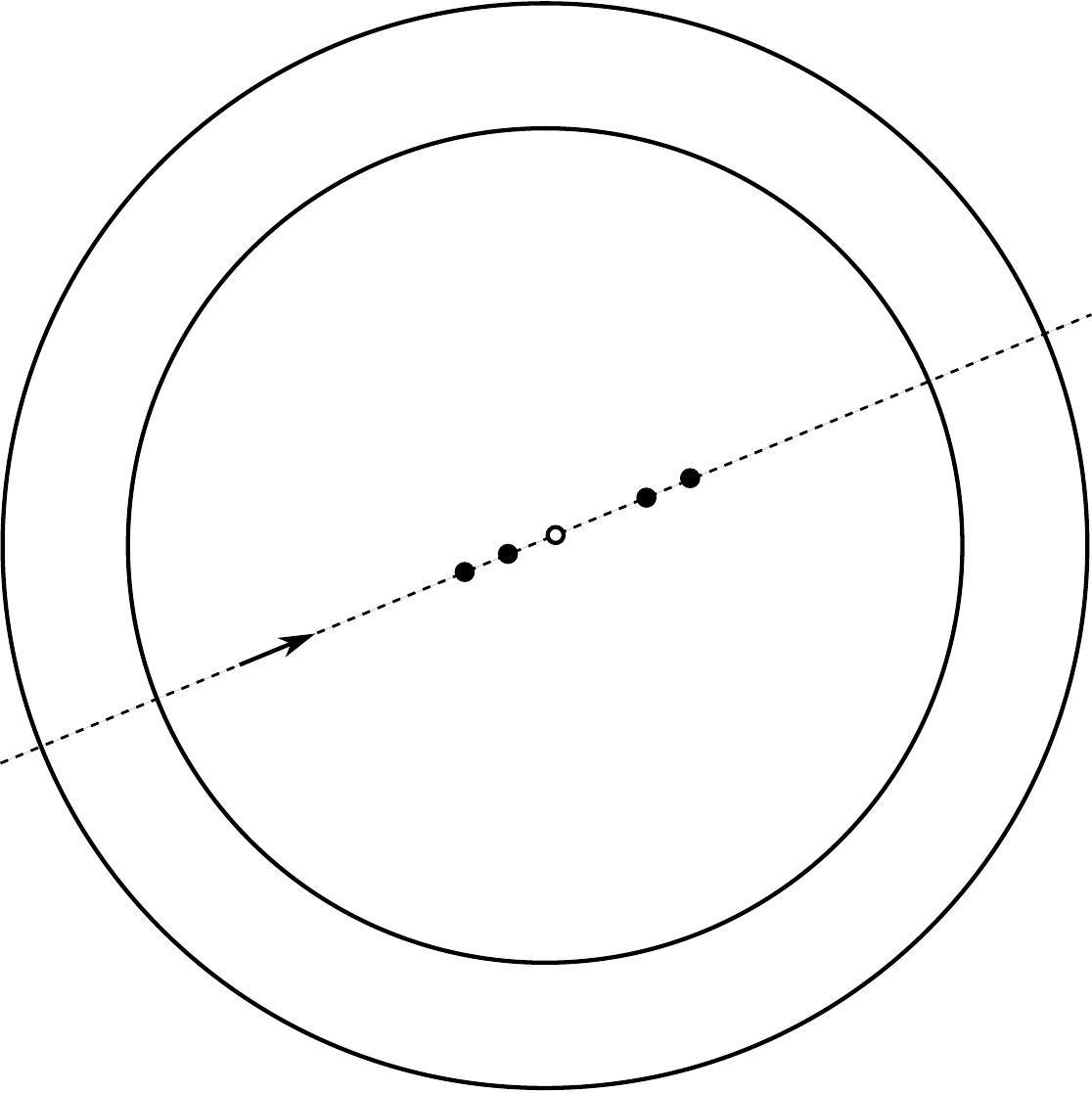}}
    \put(0.22,0.33){\color[rgb]{0,0,0}\makebox(0,0)[lb]{$\bv(l_{0,1})$}}
    \put(0,0){\includegraphics[width=\unitlength,page=2]{layer1.pdf}}
    \put(0.375,0.23){\color[rgb]{0,0,0}\makebox(0,0)[lb]{$\cgrow r$}}
    \put(0,0){\includegraphics[width=\unitlength,page=3]{layer1.pdf}}
    \put(0.351,0.069){\color[rgb]{0,0,0}\makebox(0,0)[lb]{$r^{a_1}$}}
    \put(0,0){\includegraphics[width=\unitlength,page=4]{layer1.pdf}}
    \put(0.563,0.61){\color[rgb]{0,0,0}\makebox(0,0)[lb]{$\phi$}}
    \put(0.633,0.88){\color[rgb]{0,0,0}\makebox(0,0)[lb]{$\bl_{1,k}^s$}}
    \put(0.62,0.75){\color[rgb]{0,0,0}\makebox(0,0)[lb]{$l_{1,k}$}}
    \put(0.355,0.885){\color[rgb]{0,0,0}\makebox(0,0)[lb]{$\cgrid$}}
    \put(0.48,0.9){\color[rgb]{0,0,0}\makebox(0,0)[lb]{$\cgrid r^{a_1}$}}
    \put(0.608,0.5){\color[rgb]{0,0,0}\makebox(0,0)[lb]{$\cgrid$}}
    \put(0.43,0.63){\color[rgb]{0,0,0}\makebox(0,0)[lb]{$l_{1,q}$}}
    \put(0.33,0.47){\color[rgb]{0,0,0}\makebox(0,0)[lb]{$\bl_{0,1}^m$}}
  \end{picture}%
\endgroup%
  \end{minipage} \hfill
  \begin{minipage}[c]{.45\textwidth}
    \caption{Construction of the first layer. Dashed lines represent the ``given'' $l_0$ in $\bL_0$, and two lines $l_{1,k}$ and $l_{1,q}$ from threads $k$ and $q\neq k$ in $\bL_1$, respectively. The image is a projection onto a flat surface which is parallel to $l_{1,k}$ and $l_0$. We demand that $\phi = \arcsin\ip{\bv(l_0), \bv(l_{1,k})} \in [\alpha, \beta]$. Filled dots on the picture are ``beads'' lying on $l_0, l_{1,k}$ and $l_{1,q}$; together with the    corresponding lines they form \emph{threads}.}
    \label{fig:layer1}
  \end{minipage}
\end{figure}

More precisely, given $j \in [M]$, the $j^\text{th}$ layer corresponds to a collection of truncated cylinders that are connected to the origin in $(\cgrow^j B_r) \cap \Cset$  by at most $(j-1)$ intermediate cylinders from the previous layers. More specifically, we define the set $\frL_j$ of $j^\text{th}$ layers for $j \in [M]$ recursively as follows: 
\[
  \frL_j(\bL_{j-1}) = \left\{ \bL_j := (T_{j,k})_{k=1}^{|\bL_j|} \ \Big| \ \text{each \emph{thread} $T_{j,k} = \big((\bl_{j,k}^m)_{m=1}^{N_{j,k}}, l_{j,k}\big)$ satisfies 1.-3. below} \right\}
\]
\begin{enumerate}[itemsep=1mm, topsep=0pt]
  \item The $k^{\text{th}}$ thread $T_{j,k}$ is formed by a line $l_{j,k} \in \L$ and a $\cgrid$-grid of $N_{j,k} \geq 1$ points $(\bl_{j,k}^m)_{m=1}^{N_{j,k}}$ lying on either segment of $l_{j,k} \cap A_j$.
  \item There is $q \in \bigl[ \bL_{j-1} \bigr]$, such that $\bl_{j-1,q}^s \in \Cyl{l_{j,k}} \cap l_{j-1,q}$ for some $s \in [N_{j-1,q}]$. Also, $l_{j,k} \in \eG(\bl_{j-1,q}^s; l_{j-1,q})$, where $\eG(\by;l) := \eLind{\by, \bv(l)}{\alpha, \beta} \backslash \eL \{ \bm 0, \bx \}$ for $\by \in \R^d$ and $l \in \L$, $\alpha = 2\arctan \cgrid^{-1}$ and $\beta = \arccos\cgrid^{-1}$.
  \item For every $q \in \bigl[ |\bL_j| \bigr] \backslash \{ k \}$ and $s \in [N_{j,q}]$, we have $|\bl_{j,k}^m - \bl_{j,q}^s| > \cgrid r^{a_j}$ .
\end{enumerate}
See Figure \ref{fig:layer1} for an illustration.

In order to prove \eqref{e:local}, we start by building the first layer. For the $j^\text{th}$ layer $\bL_j$, we say that it exists if all directing lines $\{l_{j,k}\}_{k=1}^{|\bL_j|}$ in its threads belong to $\omega^-$. We can show that there is a procedure allowing to construct the first layer $\bL_1$ with $\Omega(r)$ threads. The following statements hold uniformly over all $l_0\in \eL\{\bm 0\}$ and $l_x \in \eL\{\bx\}$.

\begin{lem} \label{l:first}
  There exist $f_1 > 0$ and $c > 0$, such that
  \begin{equation} \label{e:first}
    -\log \P \{ \not\exists \, \bL_1 \in \mathfrak{L_1}(\bL_0): |\bL_1| > f_1 r \} \in \Omega(r).
  \end{equation}
\end{lem}
We can also prove similar growth in the number of connections between layers $\bL_j$ and $\bL_{j+1}$ for $j \in [M-1]$.
\begin{lem} \label{l:intermediate}
  There is $\cgrow$ such that given any $f_j>0$, there exists $f_{j+1}>0$, so that for any $j^\text{th}$ layer $\bL_j$ satisfying $|\bL_j| \geq f_j r^{b_j}$,
  \begin{equation} \label{e:intermediate}
    -\log \P \bigl\{ \not \exists\, \bL_{j+1} \in \frL_{j+1}(\bL_j): |\bL_{j+1}| > f_{j+1} r^{b_{j+1}} \bigr\} \in \Omega(r^{b_{j+1}})
  \end{equation}
  for any $j \in [M-1]$.
\end{lem}
The statements in Lemmas \ref{l:first} and \ref{l:intermediate} are very similar, but proving the latter is slightly more challenging, because unlike $\bL_0$, the layers $(\bL_j)_{j=1}^{M-1}$ contain more than one thread for large $r$. Most of the peculiarities of our definition of layers come from the necessity to accommodate for that.

We will need to consider balls of the metric $\rho(\omega^-)$ around the beads of a layer $\bL_j$.
That is, we will need to consider sets of the type
\begin{equation}
\label{e:ball_layer}
  B_c^{\rho (\omega^-)}(\bL_j):=\bigcup_{k\in [|\bL_j|]}\bigcup_{m \in [N_{j,k}]} B_c^{\rho (\omega^-)}(\bl_{j,k}^m).
\end{equation}

Once we show that $\Cyl{l_0}$ is connected to sufficiently many points in $A_M$, we bridge the last layer and $\frC_0$.
\begin{lem} \label{l:last}
  Pick $\epsilon$ and $\chi$ as in Lemma \ref{l:global}. There exist $f_M>0$ and $\cpatch$ such that given any $M^\text{th}$ layer $\bL_M$ which satisfies $|\bL_M| \geq f_M r^{b_M}$, we have
  \[
    -\log \P \Big\{ \bigl| \frC_0 \cap B_{(\cpatch + \cgrow^M + 1)r}^{\rho (\omega^-)}(\bL_M) \bigr| < r^\chi \bigr\} \in \Omega(r^\chi).
  \]
\end{lem}
\section{Proofs}\label{sec:proof}
In what follows, we will use Proposition 4.1 from \cite{tyke_bro}. We use it in the following form:
\begin{equation} \label{e:p41}
  c \leq \mu \bigl( \eL \{ \bm 0 \} \cap \eL \{ \bx \} \bigr) \, |\bx|^{d-1} \leq c + c' |\bx|^{-2}, \qquad \forall \bx \in \R^d \backslash B_4.
\end{equation}
This fact has an immediate corollary which will be useful later. It can be deduced by covering $\bx + \partial B_{cR^a}$ with $cR^{a(d-1)}$ unit balls and then using \eqref{e:p41}.
\begin{cor}
  Let $R(\bx)=|\bx|$ for $\bx \in \R^d$. For any $a \in (0,1)$ and $c>0$,
  \begin{equation} \label{e:point_sphere}
    \mu \bigl( \eL\{\bm 0\} \cap \eL(\bx + B_{c R^a})\bigr)\in O(R^{(a-1)(d-1)}).
  \end{equation} 
\end{cor}
We will also use the following corollary of the Azuma-Hoeffding inequality to provide concentration bounds.
\begin{lem} \label{l:AH}
  Suppose that a filtration $\F:=(\eF_n)_{n=1}^N$ supports a sequence of random variables $(X_n)_{n=1}^N$ such that $X_n | \eF_{n-1} \sim \text{Poisson} (\mu_n)$, where $(\mu_n)_{n=1}^N \subset \R$ is such that $\mu_n \geq \mu_0 > 0$ for some $\mu_0>0$ and any $n \in [N]$. Let $I_n = \ind \{ X_n > 0 \}$. Then,
  \begin{equation} \label{e:AH}
    -\log \P \left\{ \sum_{n=1}^N I_n \leq \frac{p}{2}N \right\} > \frac{p^2}{8}N, \qquad \text{where } \ p = 1 - \exp \{-u \mu_0 \} > 0.
  \end{equation}
\end{lem}
\emph{Proof.} Note that
\[
  \E[I_n] = 1-\exp\{-u\mu_n\} \geq p, \qquad \forall n \in [N].
\]
For every $n \in [N-1]$ we have $|I_{n+1} - I_n| \leq 1$. Let $Y_n = \sum_{m=1}^n (I_m - p)$ so that $(Y_n)_{n=1}^N$ is a sub-martingale with respect to $\F$ with increments bounded by $1$ in absolute value. We can apply the Azuma-Hoeffding inequality in order to obtain:
\[
  -\log \P \left\{ \sum_{n=1}^N I_n \leq Np/2 \right\} = -\log \P \bigl\{ Y_N - Y_0 \leq -Np/2 \bigr\} \geq \frac{(Np/2)^2}{2N} =\frac{p^2}{8}N. \fqed
\]
In the remainder of this article we prove the results stated in the Section \ref{sec:res}. We start with a proof which covers Lemmas \ref{l:first} and \ref{l:intermediate} followed by a proof for Lemma \ref{l:last}. Then we prove the remaining results in reverse order finishing with the proof of the Shape Theorem.

\emph{Proof of {\bf Lemmas \ref{l:first} and \ref{l:intermediate}}}.
Order the indices $\Big\{ \bigl\{ (m,k) \bigr\}_{m=1}^{N_{j,k}} \Big\}_{k=1}^{|\bL_j|}$ in such a way that indices corresponding to earlier beads within earlier threads appear first. Specifically, define an injection $q: (m,k) \mapsto \N$ for each $m \in [N_{j,k}]$ and $k \in \bigl[ |\bL_j| \bigr]$ so that $q(m,k) < q(m',k')$ whenever $k' > k$ and $q(m,k)<q(m',k)$ if $1 \leq m < m'\leq N_{j,k}$. For simplicity, suppose that $q$ is a bijection to $[N]$ for some integer $N$. 
We have supposed that $|\bL_j| \geq f_j r^{b_j}$, so we assume that $N\asymp r^{a_j+b_j}= r^{b_{j+1}}$ which is possible by ignoring every thread after the one with index $\floo{f_j r^{b_j}}$.
Associate a filtration $\F := 
(\eF_n)_{n=1}^N$ to the threads of $\bL_j$ as follows:
\[
  \eF_n := \sigma \Big( \bigl\{ e_L: \omega \to \omega(L) \bigl| L \in \L_n \bigr\} \Big) ,\quad \L_n := \eL\{\bl_{j,k}^m: q(m,k) \in [n]\},\quad n \in [N].
\]
Next we define recursively a sequence of tuples:
\[
  \left( \eL_n, \ \ X_n, \ \ I_n, \ \ T_n=\bigl( (\bl_{j+1, n}^m)_{m=1}^{N_{j+1,n}}, \ \ l_{j+1,n}\bigr), \ \ \by_n \right)_{n=1}^N,
\]
where, for each $n \in [N]$ we write $(m,k) = q^{-1}(n)$ and define
\begin{equation} \label{eq:layer}
  \eL_n := \eG(\bl_{j,k}^m; l_{j,k}) \backslash \Big( \bigcup_{q=1}^{n-1} \eL \bigl( B_{\cgrid r^{a_{j+1}}} + \{ \bl_{j+1,q}^m \}_{m=1}^{N_{j+1,q}} \bigr) \cup \eL \{ \bl_{j,k'}^{m'}: q(m', k') > n \} \Big),
\end{equation}
$X_n := \omega(\eL_n)$ and $I_n := \ind\{X_n > 0\}$ (we clarify \eqref{eq:layer} in the Remark presented right after this proof).

$\bullet$ If $I_n = 1$, we select an arbitrary line $l_{j+1,n} \in \eL_n$ from $\omega$. Choose either one of the two line segments in $l_{j+1,n} \cap A_{j+1}$ and let $\by_n$ be its center. Put a $\cgrid$-grid $(\bl_{j+1,n}^m)_{m=1}^{N_{j+1,n}}$ on this segment and pair it with $l_{j+1,n}$ to form a thread $T_n$.

$\bullet$ If $I_n = 0$, we let $T_n$ and $\by_n$ be empty sets.

\begin{figure}[!ht]
  \centering
  \def\svgwidth{1.\textwidth}
\begingroup
  \makeatletter
  \providecommand\color[2][]{
    \errmessage{(Inkscape) Color is used for the text in Inkscape, but the 
    package 'color.sty' is not loaded}
    \renewcommand\color[2][]{}
  }
  \providecommand\transparent[1]{
    \errmessage{(Inkscape) Transparency is used (non-zero) for the text in Inkscape, but the package 'transparent.sty' is not loaded}%
    \renewcommand\transparent[1]{}
  }
  \providecommand\rotatebox[2]{#2}
  \ifx\svgwidth\undefined
    \setlength{\unitlength}{802.56544525bp}
    \ifx\svgscale\undefined
      \relax
    \else
      \setlength{\unitlength}{\unitlength * \real{\svgscale}}
    \fi
  \else
    \setlength{\unitlength}{\svgwidth}
  \fi
  \global\let\svgwidth\undefined
  \global\let\svgscale\undefined
  \makeatother
  \begin{picture}(1,0.25)
    \put(0,0){\includegraphics[width=\unitlength,page=1]{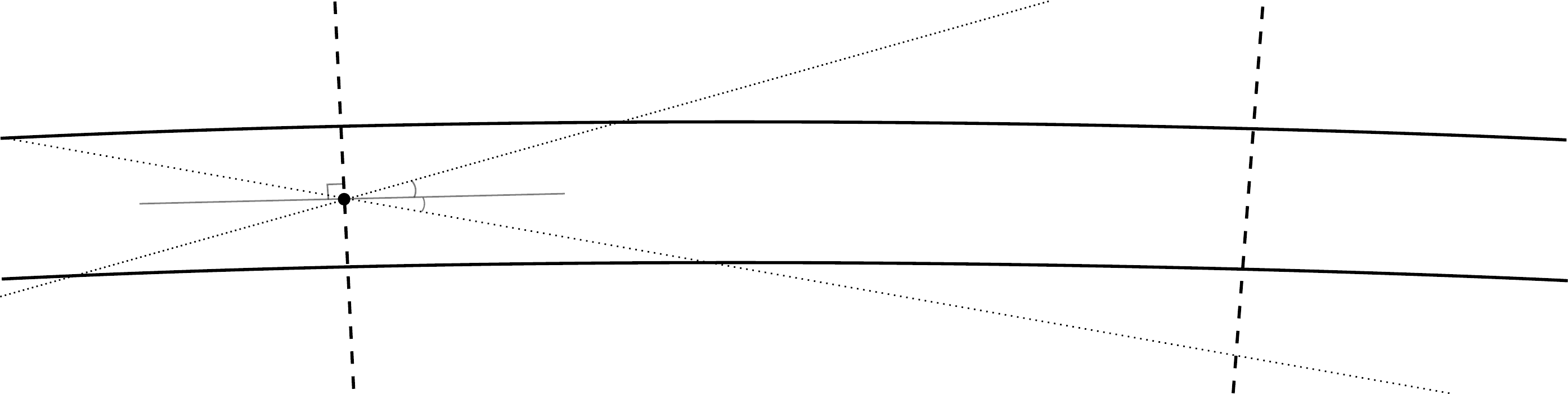}}
    \put(0,0){\includegraphics[width=\unitlength,page=2]{proof.pdf}}
    \put(0,0){\includegraphics[width=\unitlength,page=4]{proof.pdf}}
    \put(0.69,0.11){\color[rgb]{0,0,0}\makebox(0,0)[lb]{$A_j$}}
    \put(0.226,0.033){\color[rgb]{0,0,0}\makebox(0,0)[lb]{$l_{j,k}$}}
    \put(0.188,0.083){\color[rgb]{0,0,0}\makebox(0,0)[lb]{$\bl_{j,k}^m$}}
    \put(0.67,0.176){\color[rgb]{0,0,0}\makebox(0,0)[lb]{$\partial B_{\cgrow^j 
    r}$}}
    \put(0.65,0.043){\color[rgb]{0,0,0}\makebox(0,0)[lb]{$\partial B_{\cgrow^j 
    r - r^{a_j}}$}}
    \put(0.51,0.12){\color[rgb]{0,0,0}\makebox(0,0)[lb]{$r^{a_j}$}}
    \put(0.276,0.128){\color[rgb]{0,0,0}\makebox(0,0)[lb]{$\alpha$}}
    \put(0.806,0.200){\color[rgb]{0,0,0}\makebox(0,0)[lb]{$l_{j,q}$}}
    \put(0.45,0.02){\color[rgb]{0,0,0}\makebox(0,0)[lb]{at least $ \cgrid 
    r^{a_j}$}}
  \end{picture}%
\endgroup%
  \caption{For large $\cgrow$ threads become nearly orthogonal to the boundaries of $A_j$. If $l \in \eG(\bl_{j,k}^m; l_{j,k})$, then it must form an angle larger than $\alpha$ with any plane orthogonal to $l_{j,k}$. In particular, for $\cgrow$ large enough $l \cap A_j$ consists of two segments separated by distance $\Omega(r)$ from each other. The length of each segment is less then $c_0r^{a_j}$, and $\Cyl{l}$ can pass through at most one thread in $\bL_j \backslash T_k$ that is at least $\Omega(r)$ apart from $T_k$.}
  \label{fig:pf}
\end{figure}

Choose $\cgrow$ big enough so that for any $j \in [M-1]$ and $l \in \eG(\bl_{j,k}^m; l_{j,k})$, the only threads in $\bL_j$ that can be intersected by $\Cyl{l}$, except for the $k^{\text{th}}$ thread, are at a distance larger than $cr$ from $\bl_{j,k}^m$ for some $c$. See Fig.\ \ref{fig:pf} for clarification.

We can now show that $\min_{n \in [\eta N]}\mu(\eL_n) \in \Omega(1)$ for some $\eta \in (0,1)$. In fact,
\begin{equation}\label{eq:exclusion}
  \begin{split}
    \mu(\eL_n) \;&\geq \mu(\eLind{\bl_{j,k}^m, \bv(l_{j,k})}{\alpha, \beta}) - \mu(\eL \{\bl_{j,k}^m\} \cap \eL \{ \bx \}) -\mu(\eL \{\bl_{j,k}^m\} \cap \eL \{ \bm 0 \}) \\
    &-\sum_{(m',k'): q(m',k') \in (n, N]} \mu\bigl( \eLind{\bl_{j,k}^m, \bv(l_{j,k})}{\alpha, \beta} \cap \eL\{ \bl_{j,k'}^{m'}\} \bigr)\\\
    & - \sum_{q\in[n-1]: \by_q\neq \emptyset} \mu \bigl( \eLind{\bl_{j,k}^m, \bv(l_{j,k})}{\alpha, \beta} \cap \eL(\by_q + B_{\cgrid r^{a_{j+1}}}) \bigr).
  \end{split}
\end{equation}
The first term on the right-hand side is in $\Omega(1)$. The following two terms are in $O(R^{1-d})$ and $O(r^{1-d})$, respectively, due to \eqref{e:p41}. By definition, the fourth term is either zero (for $j=0$) or involves at most $N \in O(r^{b_{j+1}})$ terms all bounded above by $cr^{1-d}$ for some $c>0$, therefore the whole sum is in $O(r^{-(d-1)a_{j+1}})$. The last sum involves up to $O(r^{b_{j+1}})$ terms and each one in $O(r^{(d-1)(a_{j+1}-1)})=O(r^{-b_{j+1}})$ due to \eqref{e:point_sphere}. Therefore, there must exist an $\eta \in (0,1)$ such that $\min_{n \in [\eta N]} \mu(\eL_n) \in \Omega(1).$

We can now apply Lemma \ref{l:AH} to see that
\[
  -\log \P \left\{ \sum_{n=1}^{\floo{\eta N}} I_n \leq f_{j+1}r^{b_{j+1}} \right\} \in \Omega(r^{b_{j+1}})
\]
for some $f_{j+1}>0$. Finally, let
\[
  \bL_{j+1} = \{T_n: n \in [\eta N], \ I_n = 1\},
\]
and order it arbitrarily. $\qed$
\begin{rem}
  Roughly speaking, when $q(m,k) =n$, $\eL_n$ is a set of lines that
  \begin{itemize}[itemsep=0mm, topsep=0pt]
    \item do not hit $B_1$ or $\bx+B_1$, but go through $\bl_{j,k}^m + B_1$, and form angle $\phi \in [\alpha,  \beta]$ with a flat surface orthogonal to $l_{j,k}$;
    \item are separated by at least $c_0r^{a_{j+1}}$ from threads formed in $\bL_{j+1}$ before step $n$;
    \item stay away from the beads in $\bL_j$ that follow after step $n$.
  \end{itemize}
  When $j=0$, the third item is redundant, because $|\bL_j|=1$ and $\phi \leq \beta$.
\end{rem}

\emph{Proof of {\rm \bf Lemma \ref{l:last}}}.
Let $N = |\bL_M| \geq f_Mr^{b_M}$. Since $\bL_M$ is a layer, its threads must be separated from each other by $\cgrid r^{a_M}$, thus $N \in O(r^{d-1}) = O \bigl( |\frC_0| \bigr)$. For each thread $k \in [N]$ select only the point $\bl_k$ closest to the outer boundary of $A_M$, that is, $\partial B_{\cgrow^M r}$. Let $P := (\bl_k)_{k=1}^N$, $A := B_{\cgrow^M r} \backslash B_{\cgrow^M r - \cgrid}$ and note that $P \subset A$.  We demand that $\cpatch$ is large enough, so that for any $l \in \eL(A)$ there is at most one point in $\Cyl{l} \cap \frC_0$; for example, $\cpatch > 10\cgrow^M d$ would suffice. Associate to the points in $P$ a filtration $\F = (\eF_n)_{n=1}^N$:
\[
  \eF_n := \sigma \Big( \bigl\{ e_L: \omega \to \omega(L) \bigl| L \in \L_n \bigr\} \Big), \qquad \L_n := \eL \bigl( \{\bl_k\}_{k=1}^n \bigr), \qquad n \in [N].
\]
We build a sequence $(S_n, \eL_n, X_n, I_n, \bp_n)_{n=1}^N$ recursively. Let
\begin{eqnarray*}
  & \eL_n = \eLind{\bl_n, \bu_n}{\gamma, \pi/2} \cap \eL(\frC_0) \backslash \eL \Big( \{\bm 0, \bx\} \cup S_n \cup P \backslash \{\bl_n\} \Big), \qquad & \bu_n := \bl_n \slash |\bl_n|, \\
  & S_n := \{\bp_q: q \in [n - 1],\ I_q = 1\}, \qquad X_n = \omega(\eL_n), \qquad & I_n = \ind \{ X_n > 0\},
\end{eqnarray*}
where $\gamma > 0$ is some small angle, such that any $l \in \eLind{\bl_n, \bu_n}{\gamma, \pi/2}$ can go through at least $|\frC_0|/2$ points on $\frC_0$ for any $n \in N$.

$\bullet$ If $I_n = 1$, we take any $l \in \omega$ from $\eL_n$ and let $\bp_n = \Cyl{l} \cap \frC_0$.

$\bullet$ If $I_n = 0$, we let $\bp_n = \emptyset$.

Since $\bL_M$ is the $M^\text{th}$ layer, and since $\{\bp_n\}_{n=1}^N$ is disjoint and contained  in $\frC_0$, we have
\[
  \sum_{n=1}^N I_n \leq \Big| \frC_0 \cap  B_{(c_2 + c_4^M + 1)r}^{\rho (\omega^-)}(\bL_M) \Big|.
\]
We can show that $\mu(\eL_n)$ is bounded away from $0$:
\begin{align*}
  \mu(\eL_n) \geq \sum_{\bp \in \frC_0 \backslash S_n} & \mu \bigl( \eLind{\bl_n, \bu_n} {\gamma, \pi/2} \cap \eL \{ \bp \} \bigr) - \mu(\eL \{ \bl_n \} \cap \eL \{ \bm 0 \})  \\
  - & \mu(\eL \{ \bl_n \} \cap \eL \{ \bx \}) - \mu \bigl( \eLind{\bl_n, \bu_n} {\gamma, \pi/2} \cap \eL(\frC_0) \cap \eL( P \backslash \{\bl_n\}) \bigr).
\end{align*}
Since $N \in O(r^{d-1})$, we have $|\frC_0 \backslash S_n| \in \Omega(r^{d-1})$ and the first sum on the right-hand side is in $\Omega \bigl( |\frC_0|/r^{d-1} \bigr) = \Omega(1)$ due to \eqref{e:p41}. Inequality \eqref{e:p41} also implies that the second term is in $O(r^{1-d})$ and the third is in $O(R^{1-d})$.

We now show that the fourth term decays as $r \upinf$. Define a ``cone'' $L = \eLind{\bl_n, \bu_n}{\gamma, \pi/2} \cap \eL(\frC_0)$. Since $\gamma > 0$, any $l \in L$ intersects $A$ twice, the lengths of both segments are in $O(1)$ and the distance between them is in $\Omega(r)$. Points in $P$ are separated by at least $r^{a_M}$ with $a_M > 0$, so $l \in L$ implies that $l$ can intersect at most one point in $P \backslash \{\bl_n\}$. On the other hand, the area of $\{\by \in l \cap \partial B_{\cgrow^Mr}: l \in L\}$ is in $O(|\frC_0|^{d-1}) = O(r^{d-1})$, therefore it contains at most $O(r^{d - 1 - a_M})$ many points from $P$. The distance between $\bl_n$ and other points in $P$ that could lie in $\Cyl{l}$ for $l \in L$ is in $\Omega(r)$, therefore due to \eqref{e:p41}, we have $\mu\bigl( L \cap \eL(P \backslash \{ \bl_n\})\bigr) \in O(r^{-b_M})$.

We have thus proved that $\min_{n \in [N]} \mu(\eL_n) \in \Omega(1)$, and Lemma \ref{l:AH} now implies that
\[
  -\log \P \left\{ \sum_{n=1}^N I_n \leq cN \right\} \in \Omega(N)
\]
for some $c$. In particular, since $N \in \Omega(r^{b_M})$ and $b_M > \chi$,
\[
  -\log \P \left\{ \bigl| \frC_0 \cap B_{(c_2 + c_4^M + 1)r}^{\rho (\omega^-)}(\bL_M) \bigr| < r^\chi \right\} \geq -\log \P \Big\{ \sum_{n=1}^N I_n < r^\chi \Big\} \in   \Omega(r^\chi). \fqed
\]

\emph{Proof of {\rm \bf Proposition \ref{p:local}}}.
Split $\omega\sim\P^u$ into $M + 1$ i.i.d.\ Poisson cylinder processes of intensity $w = u/(M + 1)$, that is, we couple $\omega$ with  $(\omega_j)_{j=1}^{M+1} \iid \P^w$ such that $\sum_{j=1}^{M+1}\omega_j = \omega$. Pick $\cpatch$, $(f_j)_{j=1}^M$ and build a system of layers $(\bL_j)_{j=1}^M$ recursively as in Lemmas \ref{l:first}, \ref{l:intermediate} and \ref{l:last} with $\P = \P^w$. More specifically, define a sequence of events:
\[
  E_j = E_j(\bL_{j-1}) = \{ \exists \bL_j \in \frL_j(\bL_{j-1}): |\bL_j| > f_j r^{b_j}\}, \quad j \in [M],
\]
where for each $j \in [M]$ we pick any suitable layer $\bL_j \in \frL_j (\bL_{j-1})$ with more than $f_j r^{b_j}$ points on the event that $E_j(\bL_{j-1})$ happens. 
If $\bL_M$ is well-defined, that is, on the event $\cap_{j=1}^M E_j$, let
\[
  E_{M+1} = E_{M+1}(\bL_M) = \Big\{ \bigl|\frC_0 \cap B_{(c_2 + c_4^M + 1)r}^{\rho (\omega^-)}(\bL_M) \bigr| < r^\chi \Big\}.
\]
Let $\clocal = 2(M+1)\cgrow^M + \cpatch + 1$ and note that 
\[
  \bigcap_{j=1}^{M+1}E_{M+1} \subseteq \left\{ |\frC_0 \cap B_{\clocal r}^{\rho(\omega^- + \delta_{l_0})}| \geq r^\chi \right\}.
\] 
We finish the proof by combining Lemmas \ref{l:first}, \ref{l:intermediate} and \ref{l:last} together with the fact that $\chi>1$ to get:
\[
  -\log \P^u \bigl\{ |\frC_0 \cap B_{\clocal r}^{\rho(\omega^-+\delta_{l_0})}| < r^\chi \bigr\} \geq - \sum_{j=1}^{M+1} \log \P^w \left[ E_j^c \Big|  \bigcap_{k=1}^{j-1} E_k \right] \in \Omega(r),
\]
where an intersection over an empty set of index equals $\mathcal{M}$.
This finishes the proof. $\qed$

\emph{Proof of {\bf Lemma \ref{l:global}}}.
As $r\upinf$, the ratio of the length of the prospective highway to the size of the patches also goes to infinity, thus for $r$ big enough, if a cylinder intersects a point in $\frC_0'$ and in $\frC_x'$, those are the only two points in $\frC_0$ and $\frC_x$ that are being intersected by that cylinder. Also, if a cylinder intersects $\bm 0$ or $\bx$ and one of the points in $\frC_0'$, then there is at most one point in $\frC_x'$ that could be intersected by the same cylinder. Formally, if $l \in \eL \{ \bm 0 \} \cap \eL(\frC_0')$, $\bm a = \Cyl{l} \cap \frC_0'$ and $b_0(\bm a) = \Cyl{l} \cap \frC_x'$, then $b_0(\bm a)$ contains at most $1$ point for $r$ big enough. Similarly, if $l \in \eL \{ \bx \} \cap \eL(\frC_0')$, $\bm a = \Cyl{l} \cap \frC_0'$ and $b_x(\bm a) = \Cyl{l} \cap \frC_x'$, then $\bigl| b_x(\bm a) \bigr| \leq 1$.

Finally, we use \eqref{e:p41} to obtain the statement of this lemma:
\begin{align*}
  -\log \P \Big\{ \omega^- \bigl( \eL(\frC_0') \cap \eL(\frC_x') \bigr) = 0 \Big\} & = -\sum_{\bm{a} \in \frC_0'}\sum_{\bm{b} \in \frC_x' \backslash \{ b_0(\bm a),     b_x(\bm a) \} } \log \P \Big\{ \omega\bigl( \eL \{\bm a, \bm b \} \bigr) = 0 \Big\} \\
  \in \Omega \Big( |\frC_0'|\, (|\frC_x'|-2)\, R^{1-d} \Big) & =  \Omega(r^{2\chi}/R^{d-1}). \fqed
\end{align*}

\emph{Proof of {\bf Proposition \ref{p:bound_rho}'}}.
Pick any $l_0 \in \eL \{ \bm 0 \}$ and $l_x \in \eL \{ \bx \}$. Divide $\omega \sim \P^u$ into $3$ i.i.d.\ Poisson cylinder processes of intensity $w = u/3$. We pick some $\epsilon > 0$ so that for $\chi(\epsilon) := (d-1)(1-\epsilon)$, we have $2a\chi-(d-1) > \delta$. Note that this is possible because $\delta \in (0, 1/2)$ is less than $(2a-1)(d-1)$.

Pick $\clocal$ as in Proposition \ref{p:local} with $\P = \P^w$. Let
\[
  \frC_0' = \frC_0 \cap B_{\clocal r}^{\rho(\omega^-)},\qquad \frC_x' = \frC_x \cap  B_{\clocal r}^{\rho(\omega^-)}(\bx ).
\]
Let $F_i = \bigl\{|\frC_i'| < r^\chi \bigr\}$ for $i \in \{1, 2\}$. Thanks to Proposition \ref{p:local}, $cr^{1\land \chi} < -\log \P^w(F_i)$ for some $c$. We can similarly apply Lemma \ref{l:global}:
\[
  -\log \P^w[F|F_1^c \cap F_2^c] > c'R^{2a\chi-(d-1)}, \qquad F := \bigl\{\omega^- \bigl( \eL(\frC_0') \cap \eL(\frC_x') \bigr) = 0 \bigr\}
\]
for some $c'>0$. Both $c$ and $c'$ can be chosen independently of $l_0$ and $l_x$. If we let $\cdist = 2(\clocal + 1)$, we will have $E_{l_1,l_2} \subset F_1 \cup F_2 \cup F$ and $-\log \P^u(E_{l_1, l_2}) > c''R^\delta$ for some $c''=c''(c, c')$ which does not depend on $l_0 \in \eL \{ \bm 0 \}$ or $l_x \in \eL \{ \bx \}$. This finishes the proof. $\qed$

\emph{Proof of {\bf Proposition \ref{p:bound_rho}}}.
Take $\cdist$ as in the Proposition \ref{p:bound_rho}'. Note that given $l_0 \in \eL \{ \bm 0 \}$ and $l_x \in \eL \{ \bx \}$, $E_{l_1, l_2}$ is independent of $\omega(\eL\{\bm 0 \})$ and $\omega \bigl( \eL\{ \bx \} \bigr)$. Apply Proposition \ref{p:bound_rho}':
\begin{equation*}
  \begin{split}
    &\P \bigl[ \rho(\bm 0, \bx) > R + c_1 r\bigl| \bm 0, \bx \in \Cset \bigr] = \\
    &\P \Big[ \bigcap \bigl\{ E_{l_0, l_x} : l_0 \in \eL \{ \bm 0 \}, l_x \in \eL \{ \bx \}, \ l_0, l_x \in \omega \bigr\} \Big| \omega \bigl( \eL \{ \bm 0 \} \bigr) > 0, \omega \bigl( \eL \{ \bx \} \bigr) > 0 \Big] \leq \\
    &\sup_{l_0 \in \eL \{ \bm 0 \}, l_x \in \eL \{ \bx \} } \P(E_{l_0, l_x}) \in O \bigl( \exp\{ -R^\delta \} \bigr). 
  \end{split}
\end{equation*}
This finishes the proof of Proposition \ref{p:bound_rho}. $\qed$

\emph{Proof of the {\bf Shape theorem}}.
Pick $\cdist$ and $\delta$ as in Proposition \ref{p:bound_rho}. It suffices to prove the theorem for $c > 2\cdist$. Define the following sequence of events: 
\[
  A_n = \bigl\{ \exists \bx \in \Cset \cap B_n \cap (d^{-1/2} \, \Z^d): \rho(\bm 0, \bx) > n + \cdist n^a \bigr\},\qquad n \geq 1.
\]
The conditional probabilities of these events under $\P [\cdot | \bm 0 \in \Cset]$ can be bounded above:
\[
  \P[A_n|\bm 0 \in \Cset] \leq \sum_{\bx \in B_n \cap (d^{-1/2} \, \Z^d)} \P \bigl[ \rho(\bm 0, \bx) > n + \cdist n^a \bigl| \bm 0, \bx \in \Cset \bigr] \, \P \{ \bx \in \Cset \}.
\]
We now prove that the sum of these probabilities in $n$ is finite. Let
\[
  p_n(\by) = \P \bigl[ \rho(\bm 0, \by) > n + \cdist n^a \bigl| \bm 0, \by \in \Cset \bigr], \qquad \by \in \R^d.
\] 
We now bound $p_n$ from above uniformly on $B_n \cap (d^{-1/2} \, \Z^d)$. By Proposition \ref{p:bound_rho}, we know that $cn^\delta < -\log p_n(\by)$ for some $c$ and any $\by \in \partial B_n$. On the other hand, $p_n(\by)$ increases in $|\by|$, because to connect $\bm 0$ to a distant point by a polygonal path one must first connect it to points lying closer.

We have thus proved that $\P[A_n|\bm 0 \in \Cset] < cn^{d-1} \exp\{-n^\delta\}$ for some $c$. Thanks to the Borel-Cantelli lemma only finitely many of the events $(A_n)_{n \geq 1}$ happen $\P[\cdot|\bm 0 \in \Cset]$-almost surely. Note that if a cylinder intersects $\bx \in \R^d$, it also intersects the point closest to $\bx$ in $d^{-1/2}\, \Z^d$. That is, defining $\by(\bx):= \text{argmin}_{\bz \in d^{-1/2}\, \Z^d}|\bx - \by|$ we have:
\[
  \bx \in \Cset \iff \by = \by(\bx)  \in \Cset.
\]
Finally, since $\rho(\bx,\by) \leq d^{-1/2}$, if $\rho(\bm 0, \by) \leq n + \cdist n^a$ and $\cdist(R_0-1)^a > d^{-1/2}$, then
\[
  \rho(\bm 0, \bx) \leq n + \cdist n^a + d^{-1/2} < n + cn^a. 
\]
This finishes the proof of the Shape Theorem. $\qed$

\end{document}